\documentclass[11pt]{article}
\usepackage{amsmath}
\usepackage{amsfonts}
\usepackage{amsthm}

\begin{document}

\def\N{\mathbb{N}}
\def\F{\mathbb{F}}
\def\Z{\mathbb{Z}}
\def\R{\mathbb{R}}
\def\Q{\mathbb{Q}}
\def\H{\mathcal{H}}

\parindent= 3.em \parskip=5pt

\centerline{\bf{ON A FAMILY OF $2$-AUTOMATIC SEQUENCES}} 
\centerline{\bf {DERIVED FROM ULTIMATELY PERIODIC SEQUENCES AND}}
\centerline{\bf{GENERATING ALGEBRAIC CONTINUED FRACTIONS IN $\F_2((1/t))$}}

\centerline {(Suites $2$-automatiques \`a colone vert\'ebrale ultimement p\'eriodique)}
\vskip 0.5 cm
\centerline{\bf{by A. Lasjaunias}}
\vskip 0.5 cm
\noindent \textbf{Warning:}This note is not intended to be officially published. The matter exposed here grew from numerous exchanges during the past two months with Yining Hu (at a very and too large distance !). My aim is to report in a first draft, and in a very private and personal way, on a curious mathematical structure.   

\par We consider a large family $\Large{\mathcal{F}}$ of infinite sequences over a finite alphabet $\mathcal{A}=\left\{a_1,a_2,...,a_k\right\}$. This family includes a celebrated example of a 2-automatic sequence on the set $\left\{a,b\right\}$, called Period-doubling sequence, which has been studied in a previous article [1].

\par Each sequence ${\bf{s}}$ in $\large{\mathcal{F}}$ is built in the following way, from another sequence ${\boldsymbol{\varepsilon}}=(\varepsilon_i)_{i\geq 0}$ over $\left\{a_1,a_2,...,a_k\right\}$, this last one being ultimately periodic. 

   Starting from the empty word $W_0$, we consider the sequence of words $(W_n)_{n\geq 0}$ such that for $n\geq 0$ we have 
   $W_{n+1}=W_n,\varepsilon_n,W_n$ (note the coma is for concatenation and it will be omitted when it is suitable). Hence, we have $W_1=\varepsilon_0$, $W_2= \varepsilon_0,\varepsilon_1,\varepsilon_0$  etc...Observe by construction that $W_{n+1}$ starts by $W_n$ and therefore we may consider $W_{\infty}$ the inductive limit of these words (i.e. the word begining by $W_n$ for all $n\geq 0$). Note that, for all $n \geq 1$, $W_n$ is a palindrome, centered in $\varepsilon_n$, of length $2^n-1$.
 
  This $W_{\infty}$ represents the sequence ${\bf{s}}$, which we may denote by $\bf{s}(\boldsymbol{\varepsilon})$. Hence, we have :

   $$\textbf{s}(\boldsymbol{\varepsilon})=\varepsilon_0\varepsilon_1\varepsilon_0\varepsilon_2\varepsilon_0\varepsilon_1\varepsilon_0\varepsilon_3\varepsilon_0\varepsilon_1\varepsilon_0\varepsilon_2\varepsilon_0....$$  
   
   \par In this family $\Large{\mathcal{F}}$, each sequence $\textbf{s}=(s_i)_{i\geq 0}$ , over $\left\{a_1,a_2,...,a_k\right\}$, generates an infinite continued fraction 
   $\alpha$ in $\F_2((1/t))$, denoted $CF(\bf{s})$, by replacing the letters $a_i$ by non-constant polynomials in $\F_2[t]$ (this choice is arbitrary, hence there is a $CF(\bf{s})$ for each choice but it is considered as unique in the sequel). For basic information on continued fractions, particularly in power series fields, the reader is refered to [2]. Hence we will write :
      
      $$\alpha=CF({\bf{s}})=[s_0,s_1,....,s_n,....]$$ where the $s_i$ are the partial quotients in $\F_2[t]$. We recall that the sequence of convergents to $\alpha$ is denoted $(x_n/y_n)_{n\geq 1}$. For $n \geq 1$, we have 
      $x_n/y_n=[s_0,...,s_{n-1}]=s_0+1/(s_1+1/....)$ , hence $x_1=s_0$, $y_1=1$ and $x_2=s_0s_1+1$, $y_2=s_1$ etc...
   
   \par To be more precise about periodic sequences, we introduce the following definition.
   
     \textbf{Definition.} \emph{Let $l\geq 0$ and $d\geq 1$ be two integers. An ultimately periodic sequence $\boldsymbol{\varepsilon}$ is called of type $(l,d)$ if :
      \newline 1) $l=0$  and $\boldsymbol{\varepsilon}$ is purely periodic, with period of length $d$, this being denoted by $\boldsymbol{\varepsilon}= (\varepsilon_0,\varepsilon_1,...,\varepsilon_{d-1})^\infty$.
      \newline 2) $l>0$ and $\boldsymbol{\varepsilon}$  is ultimately periodic, with a prefix of length $l$ and a period of length $d$, this being denoted by  
         $\boldsymbol{\varepsilon}= \varepsilon_0,\varepsilon_1,...,\varepsilon_{l-1},(\varepsilon_l,\varepsilon_{l+1},...,\varepsilon_{l+d-1})^\infty $.
     \newline Here all the $\varepsilon_i$'s are in $\mathcal{A}$ (assuming that $k\geq l+d$ to allow different values to the terms of the sequence $\boldsymbol{\varepsilon}$).}
     
     \par Let us illustrate the construction of $\bf{s}(\boldsymbol{\varepsilon})$ in two basic cases :
     \newline  1) $\boldsymbol{\varepsilon}=a,b,(c)^\infty$  then $$\textbf{s}(\boldsymbol{\varepsilon})=a,b,a,c,a,b,a,c,....=(abac)^\infty$$  Note that here $\boldsymbol{\varepsilon}$ is ultimately constant (of type $(3,1)$) and $\bf{s}(\boldsymbol{\varepsilon})$ is periodic. Consequently a basic property on continued fractions shows that $\alpha=CF(\bf{s}(\boldsymbol{\varepsilon}))$ is quadratic over $\F_2(t)$.
    \newline  2) $\boldsymbol{\varepsilon}= (a,b)^\infty$    then $$\textbf{s}(\boldsymbol{\varepsilon})=abaaabababaaaba....$$
     Here $\boldsymbol{\varepsilon}$ is of type $(0,2)$ and $\bf{s}(\boldsymbol{\varepsilon})$ is the celebrated sequence mentionned above and called Period-doubling. It has been proved that $\alpha=CF(\bf{s}(\boldsymbol{\varepsilon}))$ satisfies an algebraic equation of degree 4 with coefficients in $\F_2[t]$ (see [1]).

   \par As it happens in these two simple cases, we are going to prove in the following theorem that all sequences  $\bf{s}(\boldsymbol{\varepsilon})$ in $\large{\mathcal{F}}$ generate a continued fraction $\alpha$ which is algebraic over $\F_2(t)$. During the proof, the algebraic equation satisfied by $\alpha$ will appear explicitely. 
   
     \par  Let $\boldsymbol{\varepsilon}$ be a sequence of type $(l,d)$ then we denote by $\F_2(\boldsymbol{\varepsilon})$ the subfield of $\F_2(t)$ generated by the vector $(\varepsilon_0,...,\varepsilon_{l+d-1})$ whose $l+d$ coordinates belong to $\F_2[t]$.
   
      \textbf{Theorem.} \emph{Let $l\geq 0$ and $d\geq 1$ be integers. Let $\boldsymbol{\varepsilon}$ be an ultimately periodic sequence of type $(l,d)$. Let $\bf{s}(\boldsymbol{\varepsilon})$ the sequence
      in $\F_2[t]$ and $\alpha=CF(\bf{s}(\boldsymbol{\varepsilon}))$ the continued fraction in $\F_2((1/t))$, both be defined as above. Then there is a polynomial $P$ in $\F_2(\boldsymbol{\varepsilon})[x]$ such that $\deg_x(P)=2^d$ and $P(\alpha)=0$. To be more precise, setting $\beta=1/\alpha$, there are $d+1$ elements in $\F_2(\boldsymbol{\varepsilon})$, $A$ and $B_k$ for $0\leq k \leq d-1$, such that
      $$\beta^{2^d}=A+\sum_{0\leq k \leq d-1}B_k \beta^{2^k}.$$  }

      \par Note that the case $d=1$ is trivial as we saw in case 1) above. Indeed, in that case $\boldsymbol{\varepsilon}$ is ultimately constant.  
      Hence $\textbf{s}(\boldsymbol{\varepsilon})=W,\varepsilon_l,W,\varepsilon_l,...=(W,\varepsilon_l)^\infty $ where $W$ is a finite (or empty) word and therefore $\alpha$ is quadratic. In the sequel, we may assume that $d\geq 2$.

     \par The proof of the Theorem lies on the existence of a particular subsequence of convergents to $\alpha$. These particular      
        convergents are linked to the structure of the word $\textbf{s}(\boldsymbol{\varepsilon})$. Indeed, it is natural to consider the truncation of the continued fraction $\alpha$ containing the partial quotients from $s_0$ up to $s_{2^n-2}$ for $n\geq 1$, thus corresponding to the finite word $W_n$, of length $2^n-1$, mentionned above. Hence, for $n \geq 1$, we set $u_n/v_n=[s_0,...,s_{2^n-2}]$. We have     
        $$(u_1,v_1)=(s_0,1) \quad \text{and} \quad (u_2,v_2)=(s_0s_1s_2+s_0+s_2,s_1s_2+1).$$ 
        The first step of the proof was introduced in our previous work concerning the particular case of the Period-doubling sequence. In [1,p. 4 Lemma 3.2.], using basic properties on continuants, we could prove that the pair $(u_n,v_n)$ satisfies a simple recurrence relation. 
     
        Indeed, for $n \geq 1$, we have
     
          $$(R) \qquad  u_{n+1}=\varepsilon_n u_n^2  \quad \text{and}\quad v_{n+1}=\varepsilon_n u_n v_n +1, $$  with $(u_1,v_1)=(\varepsilon_0,1)$. From $(R)$ we get immediately $$v_{n+1}/u_{n+1}=v_n/u_n + 1/u_{n+1}$$ and therefore we obtain
          
            $$ v_n/u_n=\sum_{1\leq i \leq n} 1/u_i \quad \text{ for}\quad  n\geq  1.$$ Let us consider $\beta=1/\alpha$.  Then $\beta=lim_n (v_n/u_n)$ and consequently we have 
            $$\beta=\sum_{n\geq 1} 1/u_n. \qquad  eq(0)$$
         
         \par From $eq(0)$, we will show that $\beta$ is algebraic in the following way. By successive elevation to the power 2, for $0\leq i\leq d$,  we can define inductively $d+1$ sequences, $(\varepsilon (n,i))_{n\geq 0}$, as follows:
         $$\varepsilon (n,0)=1\quad \text{and}\quad  \varepsilon (n,i+1)= \varepsilon (n,i)^2\varepsilon_{n+i}\quad \text{for}\quad 0\leq i\leq d-1.$$ 
          Note that, for $n \geq 0$, we have  $\varepsilon (n,1)= \varepsilon_n$. Then we observe that, by elevating $eq(0)$ to the power 2, we get          
           $$ \beta^2=\sum_{n\geq 1} 1/u_n^2 = \sum_{n\geq 1} \varepsilon_n / u_{n+1} = \sum_{n\geq 1} \varepsilon(n,1)/u_{n+1}. \qquad eq(1)$$
           Moreover, by successive elevation to the power 2, starting from $eq(0)$ and introducing the sequences $(\varepsilon (n,i))_{n\geq 0}$, we also get, for $0\leq i \leq d$,
         $$ \beta^{2^i}=\sum_{n\geq 1} \varepsilon(n,i)/u_{n+i}. \qquad eq(i)$$
         \textbf{ Remark.} \emph{ For all $1\leq i\leq d$ the sequence  $(\varepsilon (n,i))_{n\geq 0}$ is ultimately periodic of type $(l,d)$.
         \newline Proof by induction. This is true for $i=1$. If  $(\varepsilon (n,i))_{n\geq 0}$ is a periodic sequence of type $(l,d)$, then we have,
         for $n\geq l$,  $\varepsilon (n+d,i)= \varepsilon (n,i)$ and consequently $\varepsilon (n+d,i+1)= \varepsilon (n+d,i)^2\varepsilon_{n+d+i}= \varepsilon (n,i)^2\varepsilon_{n+i}=\varepsilon (n,i+1)$.}
         \par Now we introduce a partition of the set of positive integers into $d+1$ subsets: first the finite set $F=\{k\quad \vert \quad 1\leq k\leq l+d-1\}$ and the $d$ subsets 
         $$E_j=\{md+l+j \quad \vert \quad m\geq 1 \} \quad  \text{for} \quad 0\leq j \leq d-1.$$
         $$\N^*=\bigcup_{j=0}^{d-1} E_j \bigcup F.$$
         Linked to this partition, we introduce $d$ elements, $\beta_j$ for $0\leq j\leq d-1$, in $\F_2((1/t))$, defined by
         $$ \beta_j=\sum_{n \in E_j} 1/u_n. \qquad (B)$$
         Combining $eq(0)$ and $(B)$, and defining $z_0 \in \F_2(\boldsymbol{\varepsilon})$ by $\sum_{k \in F}1/u_k$, we can write,
         $$\beta=z_0 +\sum_{0\leq j\leq d-1} \sum_{n \in E_j} 1/u_n=z_0+ \sum_{0\leq j\leq d-1} \beta_j. \qquad Eq(0)$$
         Using the above remark, concerning the periodicity of the sequences $(\varepsilon (n,i))_{n\geq 0}$, for $0\leq j\leq d-1$ and for $1\leq i\leq d$, we can write , 
          $$\varepsilon (d+l+j-i,i)\beta_j=\sum_{m\geq 1}\varepsilon (md+l+j-i,i)/u_{md+l+j}$$
                   $$\varepsilon (d+l+j-i,i)\beta_j=\sum_{n+i\in E_j}\varepsilon (n,i)/u_{n+i}.$$
         Consequently, we observe that there exists $z_i \in \F_2(\boldsymbol{\varepsilon})$, a finite sum of the first terms in the series appearing in $eq(i)$, such that, for $1\leq i\leq d$, $eq(i)$ becomes the following equality
         $$\beta^{2^i}=z_i+ \sum_{0\leq j\leq d-1} \varepsilon (d+l+j-i,i)\beta_j. \qquad Eq(i)$$ 
         (Note that these quantities $z_i$ have a different form depending on the triplet $(l,d,i)$. See the three examples below.)
         \par Now, let us introduce the following square matrix of order $d$ :
         $$ M(d)=(m_{i,j})_{0\leq i,j\leq d-1} \quad \text{where}\quad m_{i,j}=\varepsilon (d+l+j-i,i).$$ 
         Introducing two column vectors $B$ and $C$, we observe that the $d$ equations $Eq(i)$, for $0\leq i\leq d-1$, can be summed up introducing  the following linear system $(S):\quad M(d).B=C\quad $, where   
  
          \[B = \begin{bmatrix} 
                    \beta_0 \\
                    \beta_1 \\
                    \vdots \\
                    \beta_{d-1}
                  \end{bmatrix} \quad \text{and}\quad C =
                  \begin{bmatrix}
                    \beta+z_0 \\
                    \beta^2+z_1\\
                    \vdots \\
                    \beta^{2^{d-1}}+z_{d-1}
                 \end{bmatrix}. \]
           
  We introduce the determinant, $\Delta(d)$, of the matrix $M(d)$ and also the determinant $\Delta(j,d)$ obtained from $\Delta(d)$ by replacing the column vector of rank $j$ by the column vector $C$. Hence, applying Cramer's rule for solving the linear system $(S)$, we get
  $$ \beta_j= \Delta(j,d)/\Delta(d) \quad \text{for} \quad 0\leq j\leq d-1.$$ 
  Finally, reporting these values for $\beta_j$ in $Eq(d)$, we obtain
  $$\beta^{2^d}= z_d + (\sum_{0\leq j\leq d-1}\varepsilon (l+j,d) \Delta(j,d))/\Delta(d).\quad Eq(*)$$ 
  We observe that $\Delta(d)$ belongs to $\F_2(\boldsymbol{\varepsilon})$. While, $\Delta(j,d)$ belongs to $\F_2(\boldsymbol{\varepsilon})[\beta]$. Indeed, developping the determinant $\Delta(j,d)$ along the column of rank $j$, we get $\Delta(j,d)=c_j+\sum_{0\leq k\leq d-1} b_{k,j}\beta^{2^k}$. Consequently $Eq(*)$ can be witten as expected ( with coefficients in $\F_2(\boldsymbol{\varepsilon})$) :
  $$\beta^{2^d}=A+\sum_{0\leq k \leq d-1}B_k \beta^{2^k}.\qquad Eq(**)$$ 
  So the proof of the theorem is complete.
 
 \par We present, here below, three examples. In order to avoid unnecessary complications with the subscripts, we use $(a,b,c)$ for the letters $(\varepsilon_0,\varepsilon_1,\varepsilon_2)$.
 \newline \textbf{Example 1:} Type (0,2). Period-Doubling sequence .
Let us consider the case 2, mentioned above, where $$\boldsymbol{\varepsilon}= (a,b)^\infty \quad \text{and} \quad \beta=1/CF(\textbf{s}(\boldsymbol{\varepsilon})).$$ 
We have $(\varepsilon_0,\varepsilon_1,\varepsilon_2)=(a,b,a)$ and  $(z_0,z_1,z_2)=(1/a,0,1)$.   
  \[\mathrm{ \Delta(2)=}
 \begin{vmatrix}
 1 & 1 \\
 \varepsilon_1 &\varepsilon_2
 \end{vmatrix}=a+b,\] 
 \[
    \mathrm{\Delta(0,2)=}
   \begin{vmatrix}
   \beta +z_0 & 1 \\
   \beta^2+z_1 &\varepsilon_2
 \end{vmatrix}
 \quad \textrm{ and }\quad 
    \mathrm{\Delta(1,2)=}
   \begin{vmatrix}
   1 & \beta +z_0 \\
   \varepsilon_1 &\beta^2 +z_1
 \end{vmatrix}.
 \]
 Hence we get
  $$\Delta(0,2)=\beta^2+a\beta + 1\quad \text{and}\quad  \Delta(1,2)= \beta^2 +b\beta + b/a.$$       
  Since $\varepsilon (0,2)=\varepsilon_0^2\varepsilon_1=a^2b$ and $\varepsilon (1,2)=\varepsilon_1^2\varepsilon_2=b^2a$, $Eq(*)$ becomes   
    $$(a + b)\beta^4 = a + b +(\beta^2+a\beta + 1) a^2b+ (\beta^2 +b\beta + b/a) b^2a .$$  
   From this, we get (as expected, see [1, p 2, Th 1.1] ) :
   $$\beta^4 = 1+b(a+b)+ab( a + b)\beta  +ab \beta^2. \qquad  Eq(**)$$
   
  \noindent \textbf{Example 2:} Type (1,2). Here we have : 
  $$\boldsymbol{\varepsilon}= a,(b,c)^\infty \quad \text{and} \quad \beta=1/CF(\textbf{s}(\boldsymbol{\varepsilon})).$$ 
  We have $(\varepsilon_0,\varepsilon_1,\varepsilon_2,\varepsilon_3)=(a,b,c,b)$ and  $(z_0,z_1,z_2)=(1/a+1/ba^2,1/a^2,0)$.  
  \[\mathrm{\Delta(2)=}
   \begin{vmatrix}
   1 & 1 \\
   \varepsilon_2 &\varepsilon_3
   \end{vmatrix}=b+c \] 
   \[
      \mathrm{\Delta(0,2)=}
     \begin{vmatrix}
     \beta +z_0 & 1 \\
     \beta^2+z_1 &\varepsilon_3
   \end{vmatrix}
   \quad \textrm{ and }\quad 
      \mathrm{\Delta(1,2)=}
     \begin{vmatrix}
     1 & \beta +z_0 \\
     \varepsilon_2 &\beta^2 +z_1
   \end{vmatrix}.
   \]
    Hence we get
     $$\Delta(0,2)=\beta^2+b\beta + b/a\quad \text{and}\quad  \Delta(1,2)= \beta^2 +c\beta +1/a^2+ c/a+c/ba^2.$$
     We have $\varepsilon (1,2)=\varepsilon_1^2\varepsilon_2$ and $\varepsilon (2,2)=\varepsilon_2^2\varepsilon_3$. Consequently, $Eq(*)$ becomes 
     $$( b+c)\beta^4 =(\beta^2+b\beta + b/a) b^2c+ (\beta^2 +c\beta +1/a^2+ c/a+c/ba^2) c^2b .$$  
      From this, we get :
        $$\beta^4 = bc(b+c)/a+c^2/a^2 + bc( b+c)\beta  + bc \beta^2.  \qquad  Eq(**)$$
        (Note that changing $c$ into $a$, we have $\boldsymbol{\varepsilon}= a,(b,c)^\infty =a,(b,a)^\infty =(a,b)^\infty $ and we regain the previous example and the same algebraic equation for $\beta$ as above.)  
  
  \noindent \textbf{Example 3:} Type (0,3). Here we have : 
     $$\boldsymbol{\varepsilon}= (a,b,c)^\infty \quad \text{and} \quad \beta=1/CF(\textbf{s}(\boldsymbol{\varepsilon})).$$  
     We have $(\varepsilon_0,\varepsilon_1,\varepsilon_2,\varepsilon_3,\varepsilon_4)=(a,b,c,a,b)$ and  
     $$(z_0,z_1,z_2,z_3)=(1/a+1/ba^2,1/a^2,0,1).$$  
       \[\mathrm{\Delta(3)=}
        \begin{vmatrix}
        1 & 1 & 1\\
        \varepsilon_2 &\varepsilon_3 & \varepsilon_4 \\
        \varepsilon_1^2\varepsilon_2 & \varepsilon_2^2\varepsilon_3 & \varepsilon_3^2\varepsilon_4
        \end{vmatrix}=ba^2(a+c)+ac^2(b+c)+cb^2(a+b) \] 
        
        \[\mathrm{\Delta(0,3)=}
                \begin{vmatrix}
                \beta +z_0 & 1 & 1\\
                \beta^2+z_1 &a & b \\
                \beta^4 +z_2 & c^2a & a^2b
                \end{vmatrix}=(a+b)\beta^4+(a^2b+ac^2)\beta^2+ab(a^2+c^2)\beta+\delta_0 \] 
        \[\mathrm{\Delta(1,3)=}
                        \begin{vmatrix}
                        1 & \beta +z_0 & 1 \\
                        c & \beta^2+z_1 & b \\
                        b^2c & \beta^4 +z_2 & a^2b
                        \end{vmatrix}=(a+b)\beta^4+(a^2b+ac^2)\beta^2+ab(a^2+c^2)\beta+\delta_1 \] 
                
         \[\mathrm{\Delta(2,3)=}
                                 \begin{vmatrix}
                                 1 & 1 & \beta +z_0 \\
                                 c & a & \beta^2+z_1  \\
                                 b^2c & c^2a & \beta^4 +z_2 
                                 \end{vmatrix}=(a+c)\beta^4+(ac^2+b^2c)\beta^2+ac(c^2+b^2)\beta+\delta_2 \] 
        together with
        $$\delta_0=ab(a^2+c^2)z_0+(a^2b+c^2a)z_1+(a+b)z_2$$
        $$ \delta_1=cb(a^2+b^2)+(a^2b+b^2c)z_1+(a+b)z_2$$  
        $$ \delta_2=ac(c^2+b^2)+(c^2a+b^2c)z_1+(b+c)z_2.$$  
        Here, $Eq(*)$ becomes
        $$\beta^8= z_3+ (\sum_{0\leq j\leq 2}\varepsilon (j,3)\Delta(j,3))/\Delta(3)$$ 
         and we also have
         $$\varepsilon (0,3)=a^4b^2c, \qquad \varepsilon (1,3)=b^4c^2a \quad \text{and}\quad \varepsilon (2,3)=c^4a^2b.$$
         Finally combining these values and the four values for the determinants given above, from $Eq(*)$, we get the desired outcome : 
         $$\beta^8=A+ \sum_{0\leq k\leq 2} B_k\beta^{2^k}.\qquad Eq(**)$$     
         At last, remarkably enough, we can check that the four coefficients in this last equation do not only belong to $\F_2(\boldsymbol{\varepsilon})$ but are indeed elements in $F_2[t]$ and we have
         $$A=a^3b^2c + a^2b^2c^2 + ab^3c^2 + b^4c^2 + ab^2c^3 + abc^4 + a^2bc + ab^2c + abc^2 + c^4 + 1$$
         $$B_0= a^4b2c + a^3b^2c^2 + a^2b^3c^2 + ab^4c^2 + a^2b^2c^3 + a^2bc^4$$
         $$B_1=a^3b^2c + a^2b^2c^2 + ab^3c^2 + a^2bc^3$$
         and
         $$B_2=a^2bc + ab^2c + abc^2.$$ 
         
         \par An important and last point need to be discussed. Indeed, the reader will probably ask the following question : are the sequences, belonging to the family $\Large{\mathcal{F}}$, $2$-automatic as it is indicated in the title of this note ?
         \newline There are different ways to characterize automatic sequences. A direct way is to consider the letters of the infinite world as elements in a finite field $\F_q$ of characteristic $p$. If a power series $\gamma$ in $\F_q((1/t))$ is algebraic over $\F_q(t)$, then the sequence of its coefficients is $p$-automatic (Christol's theorem).
         \par Concerning the sequences $\bf{s}(\boldsymbol{\varepsilon})$ described above, in the general case the automaticity will result from a conjecture. First we assume that the $l+d$ elements defining the sequence are in a finite field $\F_q$ of characteristic $2$ with $q=2^s \geq l+d$, and consequently we may consider the power series $\gamma \in \F_q((1/t))$ associated to this sequence. Beginning by the trivial case $d=1$, we have observed that $\bf{s}(\boldsymbol{\varepsilon})$ is utimately periodic and therefore $p$-automatic for all $p$. Note that the power series $\gamma$, associated to it, is rational and consequently it satisfies a polynomial of degree $1$ over $\F_q(t)$. We make the following conjecture :
         
         \textbf{Conjecture.} \emph{Let $l\geq 0$ and $d\geq 2$ be integers. Let $\boldsymbol{\varepsilon}$ be an ultimately periodic sequence of type $(l,d)$. Let $\bf{s}(\boldsymbol{\varepsilon})$ $=(s_n)_{n\geq 0}$ be the sequence defined above. Then there exists a finite field $\F_q$ of characteristic $2$, containing $l+d$ elements identified with the terms of this sequence so that we may consider $\gamma=\sum_{n\geq 0}s_nt^{-n}$ in $\F_q((1/t))$ and there is a polynomial $P$ in $\F_q(t)[x]$ such that $\deg_x(P)=2^{d-1}$ and $P(\gamma)=0$.}
         \par In the simpler case $(l,d)=(0,2)$ and $\boldsymbol{\varepsilon}= (a,b)^\infty$, already considered several times above, the resulting sequence $\bf{s}(\boldsymbol{\varepsilon})$ (called Period-doubling) is well known to be $2$-automatic. More precisely, the above conjecture is true. Indeed, if we identify the pair $(a,b)$ with the pair $(0,1)$ in $\F_2$ then $\gamma =\sum_{n\geq 0}s_nt^{-n} \in \F_2((1/t))$ satisfies $\gamma^2+t\gamma=t^2/(t^2+1)$.

September 2022.

\end{document}